\documentclass[12pt]{article}
\usepackage{subeqn,amsfonts}
\newcommand{\sect}[1]{\setcounter{equation}{0}\section{#1}}

\newcommand{\be}{\begin{equation}}
\newcommand{\ee}{\end{equation}}
\newcommand{\bea}{\begin{eqnarray}}
\newcommand{\eea}{\end{eqnarray}}

\newcounter{remark}

\newcommand{\II}{{\mathbb I}}
\newcommand{\ZZ}{{\mathbb Z}}
\newcommand{\eps}{{\varepsilon}}
\newcommand{\tht}{{\vartheta}}
\newcommand{\elp}{{{\cal A}_{q,p}(\widehat{sl}(2)_c)}}
\newcommand{\tr}{{\mbox{Tr}}}
\newcommand{\qpvir}{{{\cal W}_{q,p}(sl(2))}}
\newcommand{\remark}{{\sl Remark \addtocounter{remark}{1}\theremark}}

\newtheorem{thm}{Theorem}

\begin{document}
\newpage
\pagestyle{empty}
\setcounter{page}{0}
\vfill
\begin{center}

{\Large {\bf {\sf New $\qpvir$ algebras from
%
%
the elliptic algebra $\elp$}}}

\vspace{7mm}

{\large J. Avan}

\vspace{4mm}

{\em LPTHE, CNRS-URA 280, Universit\'es Paris VI/VII, France}

\vspace{7mm}

{\large L. Frappat \footnote{On leave of absence from Laboratoire de 
Physique Th\'eorique  ENSLAPP.}}

\vspace{4mm}

{\em Centre de Recherches Math\'ematiques, Universit\'e de Montr\'eal, Canada}

\vspace{7mm}

{\large M. Rossi, P. Sorba}

\vspace{4mm}

{\em Laboratoire de Physique Th\'eorique ENSLAPP
\footnote{URA 1436 du CNRS. associ\'ee \`a l'\'Ecole Normale Sup\'erieure
de Lyon et \`a l'Universit\'e de Savoie.}}
\\
{\em Annecy-le-Vieux and ENS Lyon, Lyon, France}

\end{center}

\vfill

\begin{abstract}
We construct operators $t(z)$ in the elliptic algebra $\elp$. They close an 
exchange algebra when $p^m=q^{c+2}$ for $m\in\ZZ$. In addition they commute 
when $p=q^{2k}$ for $k$ integer non-zero, and they belong to the center of 
$\elp$ when $k$ is odd. The Poisson structures obtained for $t(z)$ in these 
classical limits are identical to the $q$-deformed Virasoro Poisson algebra, 
characterizing the exchange algebras at $p \neq q^{2k}$ as new $\qpvir$ 
algebras.
\end{abstract}

\def\abstractname{R\'esum\'e}
\begin{abstract}
On construit des op\'erateurs $t(z)$ dans l'alg\`ebre elliptique $\elp$ 
formant une alg\`ebre d'\'echange quand $p^m=q^{c+2}$ o\`u $m\in\ZZ$. 
De plus, ils commutent quand $p=q^{2k}$ pour $k$ entier non nul, et ils 
appartiennent au centre de $\elp$ lorsque $k$ est impair. Les structures 
de Poisson obtenues pour $t(z)$ dans ces limites classiques sont identiques 
\`a l'alg\`ebre de Poisson Virasoro $q$-d\'eform\'ee, caract\'erisant les
structures \`a $p \neq q^{2k}$ comme de nouvelles alg\`ebres $\qpvir$.
\end{abstract}

\vfill
\vfill

\rightline{ENSLAPP-AL-649/97}
\rightline{CRM-2482}
\rightline{PAR-LPTHE 97-22}
\rightline{q-alg/9706013}
\rightline{june 1997}

\newpage
\pagestyle{plain}

\sect{Introduction}

In a first paper \cite{nous} we started the study of classical Poisson 
algebra structures obtained from the elliptic algebra $\elp$ defined by 
\cite{FIJKMY}. We showed the existence of a center at $c=-2$ generated by 
a trace formula identical to the case of trigonometric algebras \cite{RSTS}, 
leading to a set of Poisson algebra structures containing the $q$-deformed 
Virasoro algebra constructed in \cite{FR}. They were characterized by the 
particular relative position of integration contours used to define the modes
of the generating functions $t(z)$ defined as formal series.

These results naturally lead us to consider three related questions:

1 - Can one construct other commuting subalgebras in $\elp$, defined by other 
conditions on $p$, $q$, $c$ ?

2 - Which classical limits (Poisson algebra structures) does one obtain from 
such new commuting subalgebras?

3 - Can one define quantizations of these Poisson algebras, that is, closed 
algebraic structures with a supplementary parameter $\hbar$ for which one 
can define a semi-classical limit $\hbar \rightarrow 0$ such that the 
algebra structures are abelian at $\hbar =0$ and the Poisson structures 
defined in 2 are given by the leading $\hbar$ order ?

Question 3 is not reciprocal to 1 and 2 since nothing in 1 or 2 implies that
the commuting subalgebras mentioned in 1 be obtained as ``limits'' of closed 
algebras. Indeed the Poisson algebras in \cite{nous} were constructed as 
limits of algebraic structures in $\elp$ which did not close when $c \ne -2$; 
in other words the generators 
$t(z)$ at $c \ne -2$ did not close an algebra, and the quantization procedure 
in \cite{FF2} required the introduction of a new parameter as $\hbar$. 

Our starting point will be the same trace-like operator $t(z)$ which was
defined in \cite{nous}. In this paper we are going to show that:

a) if $q^{c+2}=p^m$ for any $m \in \ZZ$, the algebra $\elp$ contains a 
quadratic subalgebra generated by $t(z)$ ;

b) if in addition one has $p=q^{2k}$ for any $k \in \ZZ$, $k \ne 0$, this 
quadratic subalgebra becomes abelian; moreover if $k$ is odd, it belongs 
to the center of $\elp$.

c) we now define $\beta$ such that $p^{1-\frac{\beta}{2}}=q^{2k}$: the 
Poisson brackets structure obtained when $\beta \rightarrow 0$ is isomorphic, 
up to a factor $km$ for $k$ odd and $-km(2m-1)$ for $k$ even, to the Poisson 
brackets structure obtained in \cite{nous} when $c=-2+\beta$, $\beta 
\rightarrow 0$. In this sense the quadratic algebras at $q^{c+2}=p^m$ build 
natural quantizations of the Poisson brackets structure in \cite{nous} or 
equivalently of the $q$-deformed Virasoro algebra; we have thereby answered 
in part the question raised in \cite{FF2} concerning the construction of 
other types of $\qpvir$ - algebras (i.e. quantized $q$-deformed Virasoro 
algebras).  

\medskip

We first recall the most important notations used in \cite{FIJKMY} and the 
results obtained in \cite{nous}.

\subsection{The elliptic quantum algebra $\elp$}

The elliptic quantum algebra $\elp$ \cite{FIJKMY,JMK} was defined as follows. 
$\elp$ is an algebra of operators $L_{\eps\eps',n}$ such that 
$L_{\eps\eps',n} = 0$ if $\eps\eps' \ne (-1)^n$ ($\eps,\eps'=+$ or $-$) and 
one sets $L_{\eps\eps'}(z) = \sum_{n\in\ZZ} L_{\eps\eps',n} z^n$ (in the 
sense of formal series) which is encapsulated into a $2 \times 2$ matrix
\be
L = \left(\begin{array}{cc} L_{++} & L_{+-} \cr L_{-+} & L_{--} \cr
\end{array}\right) \,
\ee
(therefore $L_{++}(z)$ and $L_{--}(z)$ are even while $L_{+-}(z)$ and 
$L_{-+}(z)$ are odd functions of $z$).
\\
One then defines ${\cal A}_{q,p}(\widehat{gl}(2)_c)$ by imposing the following 
constraints on $L_{\eps\eps'}(z)$:
\be
R_{12}^+(z/w) \, L_1(z) \, L_2(w) = L_2(w) \, L_1(z) \, R_{12}^{+*}(z/w) \,,
\label{eq22}
\ee
where $L_1(z) \equiv L(z) \otimes \II$, $L_2(z) \equiv \II \otimes L(z)$ and
$R_{12}^+(x)$ is given by the (suitably normalized) $R$-matrix of the eight 
vertex model found by Baxter \cite{Ba}:
\be
R_{12}^+(x) = \tau(q^{1/2}x^{-1}) \frac{1}{\mu(x)} \left(\begin{array}{cccc}
a(u) & 0 & 0 & d(u) \cr 0 & b(u) & c(u) & 0 \cr 
0 & c(u) & b(u) & 0 \cr d(u) & 0 & 0 & a(u) \cr 
\end{array}\right) \label{eq23}
\ee
The functions $a(u), b(u), c(u), d(u)$ are given by 
\be
a(u) = \frac{\mbox{snh}(\lambda-u)}{\mbox{snh}(\lambda)} \,, \quad
b(u) = \frac{\mbox{snh}(u)}{\mbox{snh}(\lambda)} \,, \quad
c(u) = 1 \,, \quad
d(u) = k\,\mbox{snh}(\lambda-u)\mbox{snh}(u) \,. \label{eq24}
\ee
The function $\mbox{snh}(u)$ is defined by $\mbox{snh}(u) = -i\mbox{sn}(iu)$ 
where $\mbox{sn}(u)$ is Jacobi's elliptic function with modulus $k$.
If the elliptic integrals are denoted by $K,K'$ (let ${k'}^2=1-k^2$),
\be
K = \int_0^1 \frac{dx}{\sqrt{(1-x^2)(1-k^2x^2)}} \qquad\mbox{and}\qquad
K' = \int_0^1 \frac{dx}{\sqrt{(1-x^2)(1-{k'}^2x^2)}} \,,
\ee
the functions $a(u), b(u), c(u), d(u)$ become functions of the
variables
\be
p = \exp\Big(-\frac{\pi K'}{K}\Big) \,, \qquad
q = - \exp\Big(-\frac{\pi\lambda}{2K}\Big) \,, \qquad
x = \exp\Big(\frac{\pi u}{2K}\Big) \,. \label{eq26}
\ee
The normalization factors in (\ref{eq23}) are chosen as follows \cite{JMK}:
\bea
&& \tau(x) = x^{-1} \frac{(qx^2;q^4)_\infty \, (q^3x^{-2};q^4)_\infty}
{(qx^{-2};q^4)_\infty \, (q^3x^2;q^4)_\infty} = 
x^{-1}{\frac{\tht_{q^4}(x^2q)}{\tht_{q^4}(x^{-2}q)}} \,, \label{tau} \\
&& \frac{1}{\mu(x)} = \frac{1}{\kappa(x^2)} \frac{(p^2;p^2)_\infty}
{(p;p)_\infty^2} \frac{\tht_{p^2}(px^2) \, \tht_{p^2}(q^2)}
{\tht_{p^2}(q^2x^2)} \,, \label{mu} \\
&& \frac{1}{\kappa(x^2)} = \frac{(q^4x^{-2};p,q^4)_\infty \, 
(q^2x^2;p,q^4)_\infty \, (px^{-2};p,q^4)_\infty \, (pq^2x^2;p,q^4)_\infty}
{(q^4x^2;p,q^4)_\infty \, (q^2x^{-2};p,q^4)_\infty \, (px^2;p,q^4)_\infty 
\, (pq^2x^{-2};p,q^4)_\infty} \,, \label{kappa}
\eea
where one defines the infinite multiple products as usual by
\be
(x;p_1,\dots,p_m)_\infty = \prod_{n_i \ge 0} (1-xp_1^{n_1} \dots p_m^{n_m})
\ee
and $\tht$ is the Jacobi Theta function:
\be
\tht_a(x) = (x;a)_\infty \, (ax^{-1};a)_\infty \, (a;a)_\infty \,. \label{tht}
\ee
The $\tht$ function satisfies in particular the following identities:
\begin{subequations}
\label{theta}
\bea
&& \tht_a(ax) = \tht_a(x^{-1}) = -x^{-1} \tht_a(x) \,, \\
&& \tht_a(a^sx) = (-a^{(s-1)/2}x)^{-s} \, \tht_a(x) \,.
\eea
\end{subequations}
In eq. (\ref{eq22}) $R^{+*}_{12}$ is defined by $R^{+*}_{12}(x,q,p) \equiv 
R^+_{12}(x,q,pq^{-2c})$.
\\
The $q$-determinant $q$-$\det L(z) \equiv L_{++}(q^{-1}z) L_{--}(z) - 
L_{-+}(q^{-1}z) L_{+-}(z)$ being in the center of 
${\cal A}_{q,p}(\widehat{gl}(2)_c)$, it can be factored out, being set to the 
value $q^{c/2}$ so as to get 
\be
\elp = {\cal A}_{q,p}(\widehat{gl}(2)_c)/ 
\langle q\mbox{-}\det L - q^{c/2} \rangle \,.
\ee
>From eqs. (\ref{eq24}), (\ref{eq26})--(\ref{tht}) it follows that $\elp$ is 
well defined if $|q|$ and $|p|$ are strictly smaller than $1$, and we shall 
restrict ourselves to this sector of the parameter space.

\subsection{The center of the elliptic quantum algebra $\elp$}

For convenience, we introduce the following two matrices:
\be
L^+(z) \equiv L(q^{c/2}z) \,, \qquad
L^-(z) \equiv \sigma^1 L(-p^{1/2}z) \sigma^1 \,,
\ee
and define the operators generated by 
\be
t(z) = \tr(L(z)) = \tr\Big(L^+(q^{c/2}z) L^-(z)^{-1}\Big) \label{tz}
\ee
In ref. \cite{nous}, we proved the following results:

\begin{thm}
For all values of $p$, $q$, $c$ the operators 
$t(z),t(w)$ satisfy an exchange relation of the type
\be
t(z)t(w) = {\cal Y}(z/w)^{i_1i_2}_{j_1j_2} ~ L(w)^{j_2}_{i_2} ~ 
L(z)^{j_1}_{i_1} \,,
\ee
where the matrix ${\cal Y}(z/w)$ is given by
\be
{\cal Y}(z/w) = \bigg(\Big(R^+_{12}(w/z) \, {R^+_{12}(q^{c+2}w/z)}^{-1} \, 
{R^+_{12}(z/w)}^{-1}\Big)^{t_2} \, {R^+_{12}(q^cz/w)}^{t_2}\bigg)^{t_2} \,.
\ee
\end{thm}

\begin{thm}
When $c=-2$, the operators $t(z)$ lie in the center of 
the algebra ${\cal A}_{q,p}(\widehat{sl}(2)_{-2})$.
In particular the matrix ${\cal Y}$ is equal to the 
$4 \times 4$ unit matrix, that is $[t(z) \, , \, t(w)] = 0$.
\end{thm}

\begin{thm}
There exists a natural Poisson structure on the center
of ${\cal A}_{q,p}(\widehat{sl}(2)_{-2})$ given by
\be
\Big\{ t(z),t(w) \Big\} = -(\ln q) \left((w/z)\frac{d}{d(w/z)} 
\ln\tau(q^{1/2}w/z) - (z/w)\frac{d}{d(z/w)}\ln\tau(q^{1/2}z/w)\right) 
~ t(z)t(w) \,, \label{ps1}
\ee
and leading to a whole set of Poisson structures for the modes 
$\displaystyle t_n = \oint_C \frac{dz}{2\pi iz} \, z^{-n} \, t(z)$.
\end{thm}

\medskip

\noindent
This set of Poisson structures is parametrized by the relative positions of 
the contours $C_z$ and $C_w$ around the origin, used to extract
$\{ t_n \, , \, t_m \}$
from (\ref{ps1}). The initial Poisson bracket structure (\ref{ps1}) must be
understood in the sense of formal series for $t(z), t(w)$. Hence, depending on
the choice of relative positions of the contours for $z$ and $w$ or 
equivalently on the choice of a formal series expansion for the meromorphic 
structure functions according to whether $|z/w| \in ]q^k, q^{k-1}[$ for some 
$k$ in $\ZZ$, one gets distinct Poisson structures for the modes $t_n$ 
labeled by $k$. This fact was also hinted at in the study of the quantized 
version \cite{FF2}.

\sect{Quadratic subalgebras in $\elp$}

We now turn to the task of identifying possible closed (eventually
abelian) algebras of trace-like generators in $\elp$. Since the generic problem
is far too vast we shall simply ask the question whether the generators $t(z)$
already defined in (\ref{tz}) may close an exchange algebra. We first prove:
\begin{thm}
For any integer $m$, if $p$, $q$, $c$ are 
connected by the relation $p^m=q^{c+2}$, the generators $t(z)$ realize an 
exchange algebra with all generators $L(w)$ of $\elp$:
\be
t(z)L(w) = F\Big(m,\frac{w}{z}\Big) L(w)t(z) \label{eq31}
\ee
where
\begin{subequations}
\label{eq32}
\bea
F(m,x) &=& \prod_{s=1}^{2m} q^{-1} \,
\frac{\tht_{q^4}(x^{2}q^2p^{-s}) \, \tht_{q^4}(x^{-2}q^2p^s)}
{\tht_{q^4}(x^{-2}p^s) \, \tht_{q^4}(x^2p^{-s})} 
\quad \mbox{for $m>0$} \,, \\
F(m,x) &=& \prod_{s=0}^{2|m|-1} q \,
\frac{\tht_{q^4}(x^2p^s) \, \tht_{q^4}(x^{-2}p^{-s})}
{\tht_{q^4}(x^2q^2p^s) \, \tht_{q^4}(x^{-2}q^2p^{-s})}
\quad \mbox{for $m<0$} \,.
\eea
\end{subequations}
\end{thm}

\medskip

\noindent
{\bf Proof:}
The proof runs much along the lines of the commutativity proof in \cite{nous}.
It is easier to formulate it in terms of $L^+(w)$:
\bea
t(z) \, L_2^+(w) &=& \tr_1\Big( L_1^+(zq^{\frac{c}{2}})^{t_1}
({L}_1^-(z)^{-1})^{t_1} \Big) L_2^+(w) \nonumber \\
&=& \tr_1 \Big( L_1^+(zq^{\frac{c}{2}})^{t_1} 
({L}_1^-(z)^{-1})^{t_1} L_2^+(w) \Big)\nonumber \\
&=&
\tr_1\Big( L_1^+(zq^{\frac{c}{2}})^{t_1} 
(R_{21}^+(q^{\frac{c}{2}}w/z)^{t_1})^{-1} L_2^+(w) 
({L}_1^-(z)^{-1})^{t_1}R_{21}^{+*}(q^{-\frac{c}{2}}w/z)^{t_1} \Big) 
\,. \label{eq33} \\
&& (\mbox{by the exchange algebra for $L^-$ and $L^+$, see eq. (3.3) of
\cite{nous}}) \nonumber 
\eea
>From the exchange algebra between $L_1^+(z)$ and $L_2^-(w)$, redefining 
$z \rightarrow q^{\frac{c}{2}}z$ and using the crossing symmetry property
$(R_{21}^+(x)^{-1})^{t_1} = (R_{21}^+(xq^{-2})^{t_1})^{-1}$ one also has:
\be
L_1^+(zq^{\frac{c}{2}})^{t_1} \, 
(R_{21}^+(q^{\frac{c}{2}}q^{-c-2}w/z)^{t_1})^{-1} \, L_2^+(w) = L_2^+(w) 
\, (R_{21}^{+*}(q^{-\frac{c}{2}}w/z)^{-1})^{t_1} \, 
L_1^+(zq^{\frac{c}{2}})^{t_1} \,. \label{eq34}
\ee
In order to use (\ref{eq34}) so as to reexpress the first three factors in 
(\ref{eq33}) one needs to set $q^{-c-2}=p^{-m}$ and to use the $p$-shift 
property of $R_{21}^+$ as:
\bea
R_{21}^+(xp) &=& \Big(\tau(xq^{\frac12}) \tau(x^{-1}q^{\frac12})
\tau(xq^{\frac12}p^{\frac12}) \tau(x^{-1}q^{\frac12}p^{-\frac12})\Big)^{-1} 
R_{21}^+(x) \nonumber \\
&\equiv & F^{-1}(x)R_{21}^+(x) \,. \label{eq35}
\eea
The origin of the condition $p^m=q^{c+2}$ is precisely in that it is the only 
one allowing for a substitution of (\ref{eq34}) in (\ref{eq33}). 
\\
One then extracts from (\ref{eq33}) the prefactor generated by the use of 
(\ref{eq35}) inside (\ref{eq34}). This prefactor reads:
\be
F(m,x) \equiv \frac{R_{21}^+(xp^{-m})}{R_{21}^+(x)} = \left\{
\begin{array}{ll}
\displaystyle \prod _{s=1}^m F(xp^{-s}) & \mbox{for $m > 0$} \,, \\ \\
\displaystyle \prod _{s=0}^{|m|-1} F(xp^s)^{-1} & \mbox{for $m < 0$} \,. \\
\end{array} \right. \label{eq37}
\ee
Therefore one has:
\be
t(z)\, L_2^+(w)=F(m, q^{\frac{c}{2}}w/z)L_2^+(w) \, 
\tr_1\Big( (R_{21}^{+*}(q^{-\frac{c}{2}}w/z)^{-1})^{t_1} 
L_1^+(zq^{{\frac{c}{2}}})^{t_1} ({L}_1^-(z)^{-1})^{t_1} 
R_{21}^{+*}(q^{-{\frac{c}{2}}}w/z)^{t_1} \Big)\,,
\ee
and the two $R$-matrices cancel due to the same mechanism as in \cite{nous}:
\be
\tr_1 \Big( R_{21} Q_1 {R'}_{21} \Big) = \tr_1 \Big( Q_1 {R'_{21}}^{t_2} 
{R_{21}}^{t_2} \Big)^{t_2}\,,
\ee
if $R$ and $R^\prime$ are c-number matrices.
\\
Finally, (\ref{eq37}) can be computed using (\ref{tau}):
\[
F(m,x) = \prod_{s=1}^{2m} q^{-1} 
\frac{\tht_{q^4}(x^2q^2p^{-s}) \, \tht_{q^4}(x^{-2}q^2p^s)}
{\tht_{q^4}(x^{-2}p^s) \, \tht_{q^4}(x^2p^{-s})}
\quad \mbox{for $m>0$} \,, 
\]
and 
\[
F(m,x) = F(|m|,x^{-1}p^{\frac12})^{-1} \quad \mbox{for $m<0$} \,,
\]
which is formula 
(\ref{eq32}).
\\
Then, recalling that $L^+(w)=L(q^{\frac{c}{2}}w)$, one gets (\ref{eq31})
as stated.
\hfill \rule{5pt}{5pt}

\medskip

\remark : For $m=0$, the relation can be realized in two ways: either 
$c=-2$, which is the case studied in \cite{nous} and leads directly to a 
center $t(z)$ ($F(m,x)=1$); or $q=\exp \left ({\frac{2i\pi\ZZ}{c+2}}\right )$, 
hence $|q|=1$, which we have decided not to consider here owing to the 
potential singularities in the elliptic functions. Hence $m=0$ will be 
disregarded from now on.

\remark : Equation (\ref{eq31}) can be interpreted as meaning that 
$t(z)$ act in a uniform way as a sort of derivation on
$\elp$. We shall comment more extensively on this fact in the conclusion.

\remark: The function $F(m,x)$ is invariant under the shift 
$p\rightarrow pq^4$ due to the periodicity properties of $\tht_{q^4}$. Hence
one can restrict the parameter space of our algebras to any set 
$p\in\, ] p_0 \,,\, p_0q^4 [$.

\medskip

An immediate corollary is:

\begin{thm}
When $p^m=q^{c+2}$, $t(z)$ closes a quadratic subalgebra:
\be
t(z)t(w) = {\cal Y}_{p,q,m}\Big(\frac{w}{z}\Big) \, t(w)t(z) \label{eq39}
\ee
where
\be
{\cal Y}_{p,q,m}(x) = \left\{ \begin{array}{ll}
\displaystyle \left[ \prod_{s=1}^{2m-1} x^{2} 
\frac{\tht_{q^4}(x^{-2}p^s) \, \tht_{q^4}(x^2q^2p^s)}
{\tht_{q^4}(x^2p^s) \, \tht_{q^4}(x^{-2}q^2p^s)} \right]^2 
& \mbox{for $m>0$} \,, \\ \\
\displaystyle \left[ \prod_{s=1}^{2|m|} x^{2} 
\frac{\tht_{q^4}(x^{-2}p^s) \, \tht_{q^4}(x^2q^2p^s)}
{\tht_{q^4}(x^2p^s) \, \tht_{q^4}(x^{-2}q^2p^s)} \right]^2 
& \mbox{for $m<0$} \,. \\
\end{array} \right.
\label{eq310} 
\ee
\end{thm}

\medskip

\noindent
{\bf Proof:}
>From (\ref{eq31}) one has:
\begin{subequations}
\label{eq311}
\bea
t(z)L^+(w) &=& F\Big(m,q^{\frac{c}{2}} \frac{w}{z}\Big) L^+(w)t(z) \,, \\
t(z)(L^-(w))^{-1} &=& F^{-1}\Big(m,-p^{\frac12} \frac{w}{z}\Big) 
(L^-(w))^{-1}t(z) \,, 
\eea
\end{subequations}
hence (recalling that $t(z) = \tr(L^+(q^{c/2}z) L^-(z)^{-1})$)
\be
t(z)t(w)= \frac{\displaystyle F\Big(m,q^c \frac{w}{z}\Big)}
{\displaystyle F\Big(m,-p^{\frac12}\frac{w}{z}\Big)} \,\, t(w)t(z) \,.
\ee
The explicit expression for $F$ in (\ref{eq32}) gives the result after 
extensive use of the two ``periodicity'' properties (\ref{theta}).
\hfill \rule{5pt}{5pt}

\medskip

\remark: When $m=1$ the exchange function in (\ref{eq310}) is exactly 
the square of the exchange function in the quantization of the $q$-deformed 
Virasoro algebra proposed in \cite{FF2}, once the replacements $q^2 
\rightarrow p$, $p \rightarrow q$, $x^2 \rightarrow x$ are done. 
This is a first indication in our context that the elliptic algebra $\elp$
appears to be 
the natural setting to define quantized $q$-deformed Virasoro algebras.

\remark: As an additional connection we notice that all exchange 
functions ${\cal Y}_{p,q,m}(x)$ obey the typical identities for the 
Feigin-Frenkel function:
\be
\begin{array}{l}
\bigg. {\cal Y}(xq^2) = {\cal Y}(x) \,, \\
\bigg. {\cal Y}(xq) = {\cal Y}(x^{-1}) \,. 
\end{array}
\label{eq313}
\ee
Our exchange algebras then appear as natural generalizations of the $\qpvir$ 
algebra in \cite{FF2}. This interpretation will be reinforced by the next 
results.

\sect {Commuting subalgebras and Poisson structures}

We now show:

\begin{thm} \label{thm6}
For $p=q^{2k}$, $k\in \ZZ \backslash \{ 0 \}$, one has
\bea
F(m,x) &=& 1 \qquad \mbox{for $k$ odd} \,, \\
F(m,x) &=& q^{-2m}x^{4m}\left[ \frac{\tht _{q^4}(x^2q^2)}{\tht _{q^4}(x^2)}
\right]^{4m} \qquad \mbox{for $k$ even} \,. 
\eea
Hence when $k$ is odd $t(z)$ is in the center of the algebra $\elp$, while
when $k$ is even $t(z)$ is {\sl not} in a (hypothetical) center of $\elp$. 
However {\sl in both cases}, one has $[t(z) \, , \, t(w)] =0$.
\end{thm}

\medskip

\noindent
{\bf Proof:}
Theorem \ref{thm6} is easily proved using the explicit expression for 
$F(m,x)$ and the periodicity properties of $\tht$-functions in (\ref{theta}).
Due to Remark 3 one needs only to consider the cases $k=1$ and $k=2$.
$k=0$ is excluded since it would lead to $p=1$ and potential singularities.
\hfill \rule{5pt}{5pt}

\medskip

This now allows us to define Poisson structures on the corresponding abelian 
algebras even though $t(z)$ is not in the center of $\elp$ for $k$ even. They 
are obtained as limits of the exchange algebra (\ref{eq39}). Reciprocally
this time, since the 
initial non-abelian structure for $t(z)$ is closed, the exchange algebras 
(\ref{eq39}) are a natural quantization of the Poisson algebras which we 
obtain. This was not the case in \cite{nous}, where no intermediate closure 
condition such as $p^m=q^{c+2}$ could be obtained.

\begin{thm}
Setting $q^{2k}=p^{1-\frac{\beta }{2}}$ for any 
integer $k \ne 0$, one defines the $k$-labeled Poisson structure as:
\begin{subequations}
\label{eq42}
\bea
\Big\{ t(z) \, , \, t(w)\Big\}_k &\equiv& \lim_{\beta \rightarrow 0} \, 
\frac{1}{\beta} \, \left[ t(z)t(w)-t(w)t(z) \right] \nonumber \\
&=& 2km \, \ln q \left\{
\frac{x^2}{1-x^2} - \frac{x^{-2}}{1-x^{-2}} + \sum_{n=0}^{\infty} 
\left [ - \frac{2x^2q^{4n}}{1-x^2q^{4n}} + \frac{2x^2q^{4n+2}}{1-x^2q^{4n+2}} 
\right. \right . \nonumber \\
&& \hspace{20mm} \left. \left. + \frac{2x^{-2}q^{4n}}{1-x^{-2}q^{4n}} 
- \frac{2x^{-2}q^{4n+2}}{1-x^{-2}q^{4n+2}} \right ] \right\} 
\quad \mbox{for $k$ odd} \,, \\
&=& -2km(2m-1)\, \ln q \left\{ 
\frac{x^2}{1-x^2} - \frac{x^{-2}}{1-x^{-2}} + \sum_{n=0}^{\infty} 
\left [ - \frac{2x^2q^{4n}}{1-x^2q^{4n}} + \frac{2x^2q^{4n+2}}{1-x^2q^{4n+2}} 
\right. \right. \nonumber \\
&& \hspace{20mm} \left. \left. + \frac{2x^{-2}q^{4n}}{1-x^{-2}q^{4n}} 
- \frac{2x^{-2}q^{4n+2}}{1-x^{-2}q^{4n+2}} \right ] \right\}
\quad \mbox{for $k$ even} \,. 
\eea
\end{subequations}
\end{thm}

\medskip

\noindent
{\bf Proof:}
We note that 
\be
\Big\{ t(z) \, , \, t(w) \Big\}_k = \frac{d{\cal Y}_{p,q,m}}{d\beta}
\Big\vert_{\beta=0} \, t(z)t(w) = \frac{d\ln{\cal Y}_{p,q,m}}{d\beta}
\Big\vert_{\beta=0} \, t(z)t(w) \,,
\ee
the two equalities coming from the fact that ${\cal Y}_{p,q,m} = 1$ when
$q^{2k}=p$. 
The proof is then obvious from (\ref{eq310}) and the definition of 
$\tht$-functions (\ref{tht}) as absolutely convergent products (for $|q|<1$), 
hence as in \cite{nous}, the series in (\ref{eq42}) are convergent and define 
univocally a structure function for $t(z)$.
\hfill \rule{5pt}{5pt}

\medskip

This formula coincides exactly with the Poisson structure of the center of 
${\cal A}_{q,p}(\widehat {sl}(2)_{-2})$, provided one reabsorbs $km$ and
$-km(2m-1)$ into the definition of the classical limit as $\beta \rightarrow 
km\beta$ for $k$ odd and $\beta \rightarrow -km(2m-1)\beta$ for $k$ even.
By the same mechanism as \cite{nous}, it leads to a rich set of Poisson 
brackets for the modes of $t(z)$ defined as 
$\displaystyle t_n = \oint_C \frac{dz}{2\pi iz} \, z^{-n} \, t(z)$ 
due to the poles of the structure function at $z/w=q^n$, $n \in \ZZ$.

\section{Concluding remarks}

Theorem 7 now provides us with an immediate interpretation of the 
quadratic structures (\ref{eq39}). Since we have seen in \cite{nous} that 
the Poisson structures derived from (\ref{eq42}) contained in particular the 
$q$-deformed Virasoro algebra (up to the delicate point of the central 
extension which is not explicit in (\ref{eq42})), the quadratic algebras 
(\ref{eq39}) are inequivalent (for different values of $m$ !) quantizations 
of the $q$-deformed Virasoro algebra, globally defined on the $\ZZ$-labeled 
2-dimensional subsets of parameters defined by $p^m=q^{c+2}$. They are thus 
generalized $\qpvir$ algebras at $c=-2+m\frac{\ln p}{\ln q}$. 
Moreover the closed algebraic relation (\ref{eq31}) may, in such a frame, 
acquire a crucial importance as a $q$-deformation of the Virasoro-current 
commutations relations. This would then provide us with the full 
$q$-deformed structure required to construct possible generalizations of 
quantum Ruijsenaars-Schneider models, following the original derivation
of $q$-deformed Virasoro algebras \cite{SKAO,AKOS}.

A better understanding of the undeformed limit $q \rightarrow 1$ would help
us to clarify this interpretation if one could indeed identify the standard
Virasoro-Kac Moody structure in such a limit. The difficulty lies in the 
correct definition of this limit for the generators $L(z)$ and $t(z)$ which
should be consistent with such an interpretation.

As in \cite{FF2} the help could come from an explicit bosonization of the 
elliptic algebra as was done for $U_q(\widehat {sl}(N)_c)$ in \cite{AOS}. 
This would also
provide us with a solution of the previously mentioned central-extension
problem. At this time a bosonized version of the elliptic algebra $\elp$ is
available only at $c=1$ \cite{FSHY} using bosonized vertex operators 
constructed in \cite{AJMP}.

It is important to note also that such vertex operators were interpreted as 
$q$-analogs of primary fields for $q$-deformed $W_N$ algebras 
\cite{FJMOP,AJMP}.
Again this establishes a connection between $q$-Virasoro and $q$-$W_N$ algebras
and the elliptic algebra $\elp$ or its very recent generalization to
arbitrary $N$ ${\cal A}_{q,p}(\widehat{sl}(N)_c)$ \cite{FSHY}, but yet only
by using explicit bosonized forms. Our connection on the other hand,
where $q$-Virasoro 
algebras are directly constructed as subalgebras of $\elp$, uses only the
abstract algebraic structure of the elliptic algebra. Extensions of such a
construction may be worth looking for, providing the framework to understand
this vertex operator-primary field connection at a purely algebraic level.

\medskip

{\bf Acknowledgements}

This work was supported in part by CNRS and Foundation Angelo della Riccia. 
L.F. is indebted to Centre de Recherches Math\'ematiques of Universit\'e de 
Montr\'eal for its kind invitation and support. M.R. and J.A. 
wish to thank ENSLAPP-Lyon for its kind hospitality, and J.M. Maillet for
discussions.

\end{document}